\documentclass{article}
\usepackage[utf8]{inputenc}

\usepackage[english]{babel}

\usepackage[utf8]{inputenc}
\usepackage[T1]{fontenc}

\usepackage{amssymb}
\usepackage{stmaryrd}
\usepackage{amsmath,amsfonts}
\usepackage[tikz]{bclogo}
\usepackage[top=4.34cm]{geometry}
\usepackage{subcaption}

\usepackage{lmodern}
\usepackage{textcomp}
\usepackage{mathtools, bm}
\usepackage{amssymb, bm}
\usepackage[unq]{unique}
\usepackage{enumerate}
\usepackage{comment}
\usepackage[breaklinks]{hyperref}

\usepackage[numbers]{natbib}  

\usepackage[breaklinks]{hyperref}
\usepackage[numbers]{natbib}

\usepackage{amsthm}

\usepackage{cleveref}

\usepackage{tikz}
\usetikzlibrary {quotes} 

\newtheorem{theorem}{Theorem}[section]
\newtheorem{lemma}[theorem]{Lemma}

\newtheorem{Prop}[theorem]{Proposition}
\newtheorem{Question}[theorem]{Question}

\newtheorem{construction}[theorem]{Construction}
\newtheorem{Corollary}[theorem]{Corollary}

\theoremstyle{definition}

\newcommand{\lawrence}[1]{\textcolor{green}{\textbf{[#1]}}}
\newcommand{\julien}[1]{\textcolor{blue}{\textbf{[#1]}}}

\undef{\abs}
\DeclarePairedDelimiterX{\abs}[1]
  {\lvert}{\rvert}{\ifblank{#1}{\,\cdot\,}{#1}}
  
\newcommand{\eps}{\varepsilon}
\newcommand{\sseq}{\subseteq}

\newcommand{\EE}{\mathbb{E}}
\newcommand{\NN}{\mathbb{N}}
\newcommand{\PP}{\mathbb{P}}
\newcommand{\ZZ}{\mathbb{Z}}

\newcommand{\defined}{\mathrel{\coloneqq}}

\undef{\set}
\DeclarePairedDelimiter{\set}{\lbrace}{\rbrace}

\undef{\emptyset}
\newcommand{\emptyset}{\varnothing}

\newcommand{\inter}{\mathbin{\cap}}

\newcommand{\from}{\colon}

\title{On interval colourings of graphs}

\author{Lawrence Hollom \footnote{\href{mailto:lh569@cam.ac.uk}{lh569@cam.ac.uk}, Department of Pure Mathematics and Mathematical Statistics (DPMMS), University of Cambridge, Wilberforce Road, Cambridge, CB3 0WA, United Kingdom} \and Julien Portier\footnote{\href{mailto:jp899@cam.ac.uk}{jp899@cam.ac.uk}, Department of Pure Mathematics and Mathematical Statistics (DPMMS), University of Cambridge, Wilberforce Road, Cambridge, CB3 0WA, United Kingdom} \and Leo Versteegen \footnote{\href{mailto:lvv23@cam.ac.uk}{lvv23@cam.ac.uk}, Department of Pure Mathematics and Mathematical Statistics (DPMMS), University of Cambridge, Wilberforce Road, Cambridge, CB3 0WA, United Kingdom}}

\date{}

\begin{document}

\maketitle

\begin{abstract}
    An interval colouring of a graph $G=(V,E)$ is a proper colouring $c\from E\to \ZZ$ such that the set of colours of edges incident to any given vertex forms an interval of $\ZZ$. The interval thickness $\theta(G)$ of a graph $G$ is the smallest integer $k$ such that $G$ can be edge-partitioned into $k$ interval colourable graphs, and $\theta(n)$ is the largest interval thickness over graphs on $n$ vertices. We show that $c \frac{\log n}{\log \log n} \leq \theta(n) \leq n^{8/9+o(1)}$ for some $c>0$. In particular this answers a question by Asratian, Casselgren, and Petrosyan. \\
    In the second part of the paper, we confirm a conjecture of Axenovich that the maximum number of colours used in an interval colouring of a planar graph on $n$ vertices is at most $3n/2-2$.

\end{abstract}

\section{Introduction}
\label{sec:intro}

The notion of interval colouring of a graph was first introduced by Asratian and Kamalian \cite{asratian1994investigation}. We say that a graph $G = (V, E)$ is \emph{interval colourable} if there exists an edge-colouring $c : E \to \mathbb{Z}$ such that for every vertex $x$, the multiset $\{c(xy) : xy \in E\}$ of the colours of the edges incident to $x$ forms an interval of $\mathbb{Z}$. Such a colouring is called \emph{an interval colouring of $G$}. Note that any interval colouring is a proper colouring, i.e. a colouring where no two adjacent edges have the same colour. 

As the example of a triangle shows, not every graph is interval colourable. More generally, every graph having an edge chromatic number greater than its maximum degree is not interval colourable. Interval colourings of some special classes of graphs have been investigated over the last years, as well as some related problems, see for instance \cite{asratian2009proper,Axenovich2002,borowiecka2021consecutive,casselgren2015interval,hanson1998interval,petrosyan2014interval,pyatkin2004interval}.

The interval thickness $\theta(G)$ of a graph $G$ is the smallest integer $k$ such that $G$ can be edge-decomposed into $k$ interval colourable graphs. The interval thickness, and in general interval colourings, are of particular interest in theoretical computer science as they are related to scheduling tasks without waiting periods. For instance, if one wants to schedule a parent-teacher conference modelled by a bipartite graph $G$ such that each meeting lasts the same amount of time and there is no waiting time between the meetings for any of the parents and any of the teachers, then the interval thickness of $G$ corresponds to the minimum number of days over which the conference spans.

Some efforts have already been undertaken in order to bound $\theta(n)$, the largest interval thickness over all graphs on $n$ vertices. Asratian, Casselgren and Petrosyan \cite{asratian2022decomposing} proved that $\theta(n) \leq 2 \lceil n/5 \rceil$. 
This result was recently improved by Axenovich and Zheng \cite{axenovich2022interval} to $\theta(n) = o(n)$ using the Szemerédi regularity lemma. 
The case for the lower bound is less well studied. 
Asratian, Casselgren and Petrosyan \cite{asratian2022decomposing} asked if for every positive integer $k$ there exists a graph $G$ such that $\theta(G)=k$. 
Up until now, it was not known whether there exist graphs $G$ satisfying $\theta(G) \geq 3$. 
Our two main results are the following.

\begin{theorem}
\label{thm:LowerBoundTheta}
There exists a universal constant $C >0$ such that for every $r \in \mathbb{N}$, there exists a bipartite graph $G$ on $C(1000r)^{1000r}$ vertices satisfying $\theta(G) \geq r$. In particular there exists a universal constant $c' >0$ such that for all $n$ we have 
\begin{align*}
    \theta(n) \geq c'\frac{\log n}{\log \log n}.
\end{align*}
\end{theorem}

\begin{theorem}
\label{thm:UpperBoundTheta}
    For all $n$ we have $\theta(n) \leq n^{8/9+o(1)}$.
\end{theorem}

We prove the \Cref{thm:LowerBoundTheta} using a probabilistic construction in \Cref{thm:ProbabilisticConstruction}, and \Cref{thm:UpperBoundTheta} by finding large regular subgraphs of any sufficiently dense graph $G$ by \Cref{lem:large_factor}, via results on when a graph must contain a $k$-factor. \\
From \Cref{thm:LowerBoundTheta} we can then proceed to demonstrate that there are in fact graphs of every integer interval thickness, resolving the question of Asratian, Casselgren, and Petrosyan \cite{asratian2022decomposing}.

\begin{Corollary}
\label{thm:ThetaSurjective}
    For any integer $k$, there is some graph $G$ such that $\theta(G)=k$.
\end{Corollary}

In the second part of this paper, rather than the minimum number of interval-colourable graphs required to edge-decompose a graph, we are interested in the number of colours used in an interval-colouring. We define, for an interval-colourable graph $G$, the integer $t(G)$ to be the maximal number of colours used in any interval colouring of $G$. Asratian and Kamalian \cite{asratian1994investigation} proved that $t(G) \leq 2n$. Axenovich \cite{Axenovich2002} improved this bound over the class of planar graphs to $t(G) \leq (11/6)n$ and conjectured that this could be further improved to $t(G) \leq (3/2)n$. We confirm Axenovich's conjecture.

\begin{theorem}
\label{thm:3/2result}
Let $G$ be a planar graph on $n$ vertices that admits an interval colouring. Then $t(G)\leq (3/2)n-2$.
\end{theorem}

\section{A lower bound on $\theta(n)$}
\label{sec:unbdd-theta}

In this section we prove \Cref{thm:LowerBoundTheta} using the following probabilistic construction. 

\begin{theorem}\label{thm:ProbabilisticConstruction}
    Let $C>0$ and set $\delta=(1000r)^{-1}$, $\eps=\delta^{r+2}$ and $n=C(1000r)^{1000r}$. Let $A_1,\dotsc,A_{r-1},A_r$ and $B$ be sets of vertices with sizes $\abs{A_i}=n\delta^i$ and $\abs{B}=n$. Let $G_r$ be the random bipartite graph on parts $B$ and $A\defined A_1\cup\cdots\cup A_r$ such that every edge between $b \in B$ and $a_i \in A_i$ is present with probability $\eps\delta^{-i}$, each choice being independent of all others. Then there exists a universal threshold $c$ such that for every $C \geq c$, we have for every $r$ that $\theta(G_r) \geq r$ with probability at least $1/2$.
\end{theorem}

\Cref{thm:LowerBoundTheta} then follows as a simple corollary of \Cref{thm:ProbabilisticConstruction}. The rest of this section is devoted to proving \Cref{thm:ProbabilisticConstruction}, for which we will need the following technical lemmas. 

Our first tool will be the following Chernoff bound of Hoeffding \cite{hoeffding1963probability} (see Theorem 23.6 of Frieze's textbook \cite{frieze2016introduction} for a more modern reference). 

\begin{theorem}[Chernoff-Hoeffding inequality] 
\label{thm:chernoff-hoeffding}
    Let $S_n=X_1+\cdots+X_n$ where $X_1,\dotsc,X_n$ are independent random variables satisfying $0\leq X_i\leq 1$ for each $i$, and let $\mu=\EE S_n$.
    Then for $t\leq\mu$,
    \begin{align*}
        \PP(\abs{S_n-\mu}\geq t)\leq 2\exp\Bigl(-\frac{t^2}{2(\mu + t/3)}\Bigr).
    \end{align*}
\end{theorem}

We say that a bipartite graph $G$ of parts $C$ and $D$ is \emph{$p$-almost biregular} if for every vertex $x \in C$ we have $0.9p|D| \leq d(x) \leq 1.1p|D|$ and for every vertex $y \in D$ we have $0.9p|C| \leq d(y) \leq 1.1p|C|$. 

We say that a bipartite graph $G$ on parts $C$ and $D$ is \emph{$(\alpha, p)$-pseudorandom} if for every $U\sseq C$ and $V\sseq D$ satisfying $\abs{U}\geq\alpha \abs{C}$ and $\abs{V}\geq\alpha \abs{D}$ we have $\abs{e(U,V)-p\abs{U}\abs{V}}\leq (\abs{U}\abs{V})^{3/4+1/10}$.

\begin{lemma}
\label{prop:key-prop}
    Let $\gamma \leq 0.1$, $0 <p <1$ and $\alpha =p^2/2r^2$. Let $H'$ be a $p$-almost biregular and $(\alpha, p)$-pseudorandom bipartite graph, and let $H$ be a graph obtained from $H'$ by deleting some $\gamma$ proportion of the edges. Let the vertex parts of $H$ be $C$ and $D$. Then in any partition of $E(H)$ into $r$ parts, $E_1,\dotsc,E_r$, there is a subgraph $K\sseq H$ of diameter 4 such that $E(K)\sseq E_i$ for some $i$, and with $\abs{K\inter C}\geq \abs{C}/2r^2$.
\end{lemma}

\begin{proof}
    For any $1\leq i\leq r$, let $e_i(X, Y)$ denote the number of edges in $E_i$ between $X$ and $Y$, $\Gamma_i(x)$ denote the vertices joined to $x$ by an edge in part $E_i$ of the given edge-partition, and let $d_i(x)$ be $\abs{\Gamma_i(x)}$. 
    Then, letting $k$ be such that $E_k$ has the most edges of any of $E_1,\dotsc,E_r$, we have
    \begin{align*}
        \sum_{x\in C}e_k(C,\Gamma_k(x))&=\sum_{y\in D}e_k(C,y)^2\\
        &\geq\abs{D}^{-1}\Bigl(\sum_{y\in D}e_k(C,y)\Bigr)^2.
    \end{align*}

    Let $m$ be the number of edges between $C$ and $D$ in $H$. 
    Then
    \begin{align*}
        \sum_{y\in D}e_k(C,y)\geq \frac{e(C,D)}{r}=\frac{m(1-\gamma)}{r}.
    \end{align*}

    Therefore there is some $x_0\in C$ such that $e_k(C,\Gamma_k(x_0))$ is at least the average. Thus
    \begin{align}
    \label{eq:lowerboundEk}
        e_k(C,\Gamma_k(x_0))\geq\frac{m^2(1-\gamma)^2}{r^2\abs{C}\abs{D}}\geq \frac{0.9^2p^2\abs{C}\abs{D}(1-\gamma)^2}{r^2}.
    \end{align}

    The final inequality above comes from $H'$ being $p$-almost biregular. 
    Now define $E\defined \Gamma_k(\Gamma_k(x_0))$.
    Noting that $E\sseq C$, it suffices to prove that $\abs{E}\geq \abs{C}/2r^2$. As $H'$ is $p$-almost biregular, we have $\abs{\Gamma_k(x_0)}=d_k(x_0)\leq d(x_0)\leq 1.1p\abs{D}$. Hence by \eqref{eq:lowerboundEk} we already have $\abs{E} \geq \frac{e_k(C,\Gamma_k(x_0))}{\abs{\Gamma_k(x_0)}} \geq \alpha \abs{C}$. Moreover, as $H'$ is $(\alpha,p)$-pseudorandom, the number of edges in $H'$ between $E$ and $\Gamma_k(x_0)$ is at most $1.1p\abs{E}\abs{\Gamma_k(x_0)}$. Therefore
    \begin{align*}
        \abs{E}\geq \frac{e_k(E,\Gamma_k(x_0))}{1.1p\abs{\Gamma_k(x_0)}} \geq \frac{e_k(C,\Gamma_k(x_0))}{1.1^2 p^2 \abs{D}} \geq \frac{0.9^2 p^2\abs{C}\abs{D}(1-\gamma)^2}{1.1^2p^2\abs{D}r^2} \geq\frac{\abs{C}}{2r^2},
    \end{align*}
    as required.
\end{proof}

We now prove the following simple lemma, which is the tool we shall use to deduce that certain graphs cannot be interval coloured.

\begin{lemma}
\label{lem:interval-colour-diameter}
    Let $G$ be a graph and let $H$ be a subgraph of $G$. Let $d$ be the diameter of $H$ and suppose that for every vertex $v$ in $H$, we have $d_G(v) \leq \Delta$ for some $\Delta > 0$. Then if there exists a graph $K$ such that $K$ admits an interval colouring $c$ and such that $H \sseq K \sseq G$, then for every pair of edges $e,e'$ in $H$, we have $\abs{c(e)-c(e')} \leq (d+1)(\Delta-1)$.

    In particular, no vertex of $K$ has more than $(d+1)(\Delta-1)$ neighbours in $H$.
\end{lemma}
\begin{proof}
    For every pair of edges $e,e'$, consider a shortest path from $e$ to $e'$ in $H$. This path has length at most $d+1$, and colours of consecutive edges along this path can differ by at most $\Delta-1$. Therefore $\abs{c(e)-c(e')} \leq (d+1)(\Delta-1)$, as required. 
\end{proof}

\begin{proof}[Proof of \Cref{thm:LowerBoundTheta}.]
Set $\alpha_i=(\eps \delta^{-i})^2/(2r^2)^i$ and $\beta=0.1$. Let $I_i$ be the event that the graph $H_i=G[B \cup A_i]$ is $(\alpha_i,\eps\delta^{-i})$-pseudorandom, let $J_i$ be the event that $H_i$ is $\eps \delta^{-i}$-almost biregular, and let $L= I_1 \cap \dots \cap I_r \cap J_1 \cap \dots \cap J_r$. 
The proof will follow in 2 stages. First we show that $L$ occurs with probability at least $1/2$. 
Then, conditionally on the event $L$, we will inductively apply \Cref{prop:key-prop} to show that $\theta(G) \geq r$. \\

Step 1: Let $U_i \sseq A_i$ and $V_i \sseq B$ satisfying $\abs{U_i} \geq \alpha_i \abs{A_i}$ and $\abs{V_i} \geq \alpha_i \abs{B}$. Then by Chernoff-Hoeffding \ref{thm:chernoff-hoeffding} we have
\begin{align*}
        \mathbb{P}\bigl(\abs{e(U_i,V_i)-\eps \delta^{-i}\abs{U_i}\abs{V_i}}\leq (\abs{U_i}\abs{V_i})^{3/4+\beta}\bigr) \leq \exp\Bigl(-\frac{(\alpha_i^2 |A_i||B|)^{1/2+2\beta}}{4\eps \delta^{-i}}\Bigr).
\end{align*}

By union bound, it follows that
\begin{align*}
        \mathbb{P}(\Bar{I_i}) &\leq |A_i|^{|A_i|}|B|^{|B|}\exp\Bigl(-\frac{(\alpha_i^2 |A_i||B|)^{1/2+2\beta}}{4\eps \delta^{-i}}\Bigr) \\
        &\leq \exp\Bigl(2\log(n)n-\frac{n^{1+1/100}}{8}\Bigr) \\
\end{align*}

For a vertex $v \in H_i$, let $d_i(v)$ be the degree of $v$ in $H_i$. Let $x \in A_i$. Then by Chernoff-Hoeffding \ref{thm:chernoff-hoeffding} we have 

$$\mathbb{P}(|d_i(x) - \eps \delta^{-i}|B|| \geq 0.1\eps \delta^{-i}|B|) \leq \exp( - \eps \delta^{-i}|B|/2 )$$

Similarly, for $y \in B$, we have:

$$\mathbb{P}(|d_i(y) - \eps \delta^{-i}|A_i|| \geq 0.1\eps \delta^{-i}|A_i|) \leq \exp( - \eps \delta^{-i}|A_i|/2 )$$

Hence by union bound, noting that $\abs{B}>\abs{A_i}$,
\begin{align*}
    \mathbb{P}(\Bar{J_i}) \leq 2|B|\exp( - \eps \delta^{-i}|A_i|/2 ) = 2|B|\exp( - \eps |B| /2 ) \leq 2n\exp(-n^{1/2}/2).
\end{align*}

By union bound once again
\begin{align*}
        \mathbb{P}(\Bar{L}) &\leq r\exp\Bigl(2\log(n)n-\frac{n^{1+1/100}}{8}\Bigr)+2rn\exp(-n^{1/2}/2) \\
        &\leq \exp\Bigl(3\log(n)n-\frac{n^{1+1/100}}{8}\Bigr)+\exp(4\log(n)-n^{1/2}/2).
\end{align*}

It is easy to see that there exists $n_0$ such that for every $n \geq n_0$, the above expression is at most $1/2$. Setting $c=n_0$, we then get $\mathbb{P}(L) \geq 1/2$. \\

Step 2: For the rest of this proof, we assume that $L$ is realised, and we suppose for contradiction that $\theta(G) \leq r-1$ and therefore fix a partition of $E(G)$ into parts $E_1,\dotsc,E_{r-1}$, each of which is interval colourable.
Using \Cref{prop:key-prop}, we will inductively construct some subsets $B_r \sseq B_{r-1}\sseq B_{r-2}\sseq\cdots\sseq B_1 \sseq B_0=B$ satisfying $\abs{B_i}\geq\abs{B_{i-1}}/2(r-i+2)^2$, and subgraphs $K_i$ of $H_i$ for each $i$.
These will be such that for each $i$, $K_i$ is a diameter-4 bipartite subgraph of $H_i$ with $K_i\cap B=B_i$, and $E(K_i)\sseq E_{f_i}$ for some $f_i\notin\set{f_1,\dotsc,f_{i-1}}$. 
As $f_1\dotsc,f_r$ are distinct integers satisfying $1\leq f_i\leq r-1$, this will contradict that $\theta(G) \leq r-1$, as desired. 

We first construct $B_1$. By \Cref{prop:key-prop} there exists a subgraph $K_1$ of $G$ with diameter 4 satisfying $\abs{K_1\cap B}\geq \abs{B}/2r^2$, and $E(K_1)\sseq E_{f_1}$ for some $f_1$.
Define $B_1\defined K_1\cap B$. 
Now note that the degree (in $G$) of each vertex of $K_1$ is bounded above by $\max(1.1r\eps n, 2\eps\delta^{-1}n)=2\eps\delta^{-1}n$. 
Moreover, note that as $K_1$ has diameter 4, we may apply \Cref{lem:interval-colour-diameter} to see that no other vertex in $G$ may have more than $8\eps\delta^{-1}n$ edges from $E_{f_1}$ into $K_1$ while preserving interval colourability of $G$ restricted to $E_{f_1}$.

Suppose now that we have constructed $B_{k-1}\sseq B_{k-2}\sseq\cdots\sseq B_1$ with $\abs{B_{i+1}}\geq\abs{B_{i}}/2(r-i+1)^2$ and subgraphs $K_i$ of $H_i$ for $1\leq i \leq k-1$ which satisfy $K_i\cap B = B_i$ and $E(K_i)\cap E=E_{f_i}$ for some distinct $f_1,\dotsc,f_{k-1}$.
Following similar reasoning to the case for constructing $B_1$, the maximum degree of $K_i$ is at most $2\eps\delta^{-i}n$, and so again by \Cref{lem:interval-colour-diameter}, no vertex of $A_k$ can have more than $8\eps\delta^{-i}n$ edges from $E_{f_i}$ into $B_i$.
Therefore each vertex of $A_k$ must have at most 
\begin{align*}
    8\eps n(\delta^{-1}+\dotsb+\delta^{1-k})=8\eps n\frac{\delta^{1-k}-1}{1-\delta}
\end{align*} 
edges from $E_{f_1}\cup\cdots\cup E_{f_{k-1}}$ into $B_{k-1}$.
Consider removing these edges.
Each vertex of $A_k$ has degree at least $\eps\delta^{-k}n/2$, and so the proportion of edges removed from each vertex of $A_k$ is then at most 
\begin{align*}
    16\frac{\delta-\delta^k}{1-\delta}\leq 17\delta.
\end{align*}
The graph remaining after deleting these edges is coloured with at most $r-k+1$ colours, and so we may apply \Cref{prop:key-prop} to find a diameter-4 bipartite subgraph $K_k$ of $H_k$ satisfying $E(K_k)\sseq E_{f_k}$ for some $f_k\notin\set{f_1,\dotsc,f_{k-1}}$. Then $B_k\defined K_k\cap B_{k-1}$ satisfies $\abs{B_k}\geq\abs{B_{k-1}}/2(r-k+2)^2$, as required.

This completes the induction, and finishes the proof.
\end{proof}

We may now use this result to deduce \Cref{thm:ThetaSurjective}.

\begin{proof}[Proof of \Cref{thm:ThetaSurjective}]
    We want to construct a graph of interval thickness $k$.
    By \Cref{thm:LowerBoundTheta} we know that there is some graph $G_0$ on $n$ vertices with $\theta(G_0)\geq k$.
    Arbitrarily order the vertices of $G_0$ as $v_1,v_2,\dotsc,v_n$, and define, for each $i$, $G_i\defined G_0[V(G_0)\setminus\set{v_1,\dotsc,v_i}]$. Then we may note that $G_i$ is the union of $G_{i+1}$ and a star graph, which is interval colourable.
    Therefore $\theta(G_{i+1})\geq \theta(G_i)-1$.
    Furthermore, if $l$ is maximal such that $G_l$ has any edges, then $G_l$ is a star, and so $\theta(G_l)=1$.
    Thus for some $0\leq i \leq l$, we have $\theta(G_i)=k$, as required.
\end{proof}

\section{An upper bound on $\theta(n)$}
\label{sec:upper-bound}

We will prove \Cref{thm:UpperBoundTheta} by reducing it to the case of bipartite graphs with parts of equal sizes.
\begin{theorem}
\label{thm:UpperBoundThetaBipartite}
    Let $G$ be a bipartite graph on $n$ vertices and equal parts. Then $\theta(G) \leq n^{8/9+o(1)}$.
\end{theorem}

\begin{proof}[Proof of \Cref{thm:UpperBoundTheta} assuming \Cref{thm:UpperBoundThetaBipartite}.]
    Take an arbitrary graph $G$ on $n$ vertices, and take integer $s$ minimal so that $n\leq 2^s$.
    Note that we may assume that $n=2^s$ by adding extra vertices of degree zero to $G$.
    We will edge-partition $G$ into $s$ bipartite graphs $G_0,\dotsc,G_{s-1}$, each with equally-sized parts, as in the statement of \Cref{thm:UpperBoundThetaBipartite}.

    We now define our partition.
    Label the vertices of $G$ as $\set{0,1,\dotsc,2^s-1}$.
    Each $G_i$ will have the same vertex set as $G$.
    Then, for any edge $uv\in E(G)$, let $i$ be the minimal integer such that $u$ and $v$ differ modulo $2^i$, so $0\leq i \leq s-1$.
    Equivalently, $i$ is the location of the least significant digit in which the binary representation of $u$ and $v$ differ.
    Assign edge $uv$ to the graph $G_i$.
    
    Note that every graph $G_i$ is bipartite with equally-sized parts; membership of the two parts of $G_i$ is decided by whether a vertex has a 0 or 1 as the $i^\text{th}$ digit of its binary representation. 
    Note further that every edge $uv$ must be put into exactly one of the $G_i$ as there is a unique first position where the binary representations of $u$ and $v$ differ.

    Therefore \Cref{thm:UpperBoundThetaBipartite} applies to each $G_i$, and thus $\theta(G)\leq \theta(G_1)+\cdots+\theta(G_s)\leq sn^{8/9+o(1)}$.
    Noting that $s<1+\log_2(n)\leq 2\log(n)$, we may immediately deduce the conclusion of \Cref{thm:UpperBoundTheta}.
\end{proof}

The rest of this section is then devoted to proving \Cref{thm:UpperBoundThetaBipartite}. We will need the following criterion for finding a $k$-factor in a bipartite graph with parts of equal size, which can be found for instance in \cite{plummer1986matching}, page 70, Thm. 2.4.2.

\begin{theorem}
\label{thm:criterion_k-factor}
    Let $G$ be a bipartite graph on $2n$ vertices in the classes $V_1$ and $V_2$, where $\abs{V_1}=\abs{V_2}=n$. Then $G$ has a $k$-factor if and only if for all $X \sseq V_1$ and $Y \sseq V_2$ we have
    \begin{align*}
    kn+e(X,Y) \geq k\abs{X}+k\abs{Y}.
\end{align*}
\end{theorem}

The arboricity of a graph $G$ is the minimum number of forests edge-decomposing $G$. Dean, Hutchinson and Scheinerman proved \cite{dean1991thickness} that the arboricity of any graph on $m$ edges is at most $\sqrt{m/2}$. As noted in \cite{asratian2022decomposing}, since forests are interval colourable, this implies the following.

\begin{lemma}
\label{lem:BoundArboricity}
    Let $G$ be a graph on $m$ edges. Then $\theta(G) \leq \sqrt{m/2}$.
\end{lemma}

We will prove the following lemma.

\begin{lemma}
\label{lem:large_factor}
    Let $\delta=1/4$. Let $G$ be a bipartite graph on $n$ vertices of density $d$ on parts of equal sizes $X$ and $Y$. Then $G$ contains a $k$-regular subgraph on at least $d^{1/\delta}n^2/200$ edges.
\end{lemma}

We now show that \Cref{thm:UpperBoundTheta} follows from \Cref{lem:large_factor}. 

\begin{proof}[Proof of \Cref{thm:UpperBoundThetaBipartite} assuming \Cref{lem:large_factor}.]
    Take $G$ bipartite on $N=2n$ vertices with equal parts, let $d$ be the density of $G$, and let $\gamma=(1/2+d/2)^{-1}$.
    We will construct a sequence of graphs $G=G_1\supseteq G_2\supseteq\dots\supseteq G_r$ of decreasing density.
    Let the density of graph $G_i$ be $d_i$.

    If $d_i\geq N^{-\gamma}$, then by \Cref{lem:large_factor} we can find a $k$-regular subgraph $H_i$ of $G_i$ with at least $d_i^{1/\delta}N^2/200 \geq N^{2-\gamma/\delta}/200$ edges. 
    We then define $G_{i+1}$ to be $G_i$ with all edges of $H_i$ removed.

    On the other hand, if $d_i< N^{-\gamma}$ then we stop the process.
    Note that, as each step removes at least $N^{2-\gamma/\delta}/200$ edges, this process must eventually terminate.
    Say that the process stops after constructing the graph $G_r$.

    We have thus constructed an edge-partition of $G$ into $H_1,H_2,\dotsc,H_{r-1}$ and $G_r$.
    Note that as each $H_i$ is regular and bipartite, we may apply Hall's theorem to see that $H_i$ is interval colourable.
    Moreover, each $H_i$ has at least $N^{2-\gamma/\delta}/20$ edges, and thus $rN^{2-\gamma/(2\delta)}\leq 200N^2$, giving $r \leq 200N^{\gamma/(2\delta)}$.
    Finally, $G_r$ has density at most $N^{-\gamma}$, so by \Cref{lem:BoundArboricity}, $\theta(G_r) \leq N^{1-\gamma/2}$. 
    Putting this all together, we see that $\theta(G) \leq r+N^{1-\gamma/2} \leq 20N^{\gamma/(2\delta)}+N^{1-\gamma/2}=201N^{1/(\delta/2+1)}=N^{8/9+o(1)}$, as required.
\end{proof}

Our main technical lemma is the following density increment. For two vertex sets $A,B$ in a graph $G$ we write $d_{A,B}$ for the edge density $\frac{e(A,B)}{\abs{A}\abs{B}}$.

\begin{lemma}
\label{lem:DensityIncrement}
    Let $G$ be a bipartite graph on vertex classes $X$ and $Y$ of size $n$ with edge denstiy $d=d_{X,Y}$. At least one of the following is true:
    \begin{itemize}
        \item There exists $X' \subset X$ and $Y' \subset Y$ such that $\abs{X}=\abs{Y'}$ and $\frac{d_{X',Y'}}{d} \geq (\frac{n}{\abs{Y'}})^{\delta}$,
        \item For every $A \sseq X$ and $B \sseq Y$ we have $kn+e(A,B) \geq k\abs{A}+k\abs{B}$ where $k=dn/100$.
    \end{itemize}
\end{lemma}

\begin{proof}
    Suppose that the statement was not true. In particular, there exist $A \sseq X$ and $B \sseq Y$ such that $kn+e(A,B) < k\abs{A}+k\abs{B}$. By symmetry, we may assume that $\abs{A} \leq \abs{B}$. Let $C= Y \setminus B$ and $D = X \setminus A$. Since we must have $\abs{A}+\abs{B}>n$, we may infer that $\abs{A}\geq \abs{C}$ and $\abs{B}\geq \abs{D}$. Since $\abs{C}+\abs{B}=n$, we also have $\abs{D}\geq \abs{C}$.

    It follows by exclusion of the first case in the claim that there can be no sets $X'\subset X$, $Y'\subset Y$ such that $\abs{Y'}\leq \abs{X'}$ and $\frac{d_{X',Y'}}{d} > (\frac{n}{\abs{Y'}})^{\delta}$. Indeed, if this were the case, we could simply sample a set $X''\sseq X'$ of size $Y'$ at random such that $\frac{d_{X'',Y'}}{d}\geq \frac{d_{X',Y'}}{d} > (\frac{n}{\abs{Y'}})^{\delta}$.

    Therefore, we must have $d_{A,C}/d<(n/\abs{C})^\delta$ and $d_{D,Y}<(n/\abs{D})^\delta$, or using numbers of edges rather than densities, $e(A,C)<d\abs{A}\abs{C}^{1-\delta}n^{\delta}$ and $e(D,Y)<d\abs{D}^{1-\delta}n^{1+\delta}$. Combining this with $e(A,B) < k\abs{A}-k\abs{C}$, we obtain

    \begin{align*}
        e(X,Y)&=e(A,B)+e(A,C)+e(D,Y)\\
        &< k(\abs{A}-\abs{C})+d\abs{A}\abs{C}^{1-\delta}n^{\delta}+d\abs{D}^{1-\delta}n^{1+\delta}\\
        &=dn^2\left(\frac{\abs{A}-\abs{C}}{100n}+\left(\frac{n-\abs{A}}{n}\right)^{1-\delta} +\frac{\abs{A}}{n}\left(\frac{\abs{C}}{n}\right)^{1-\delta}\right).
    \end{align*}

    Letting $x=\abs{A}/n$ and $y=\abs{C}/n$, the above becomes

    \begin{align*}
        e(X,Y)<dn^2\left(\frac{x-y}{100}+({1-x})^{1-\delta} +xy^{1-\delta}\right)
    \end{align*}
    
    However, noting that $y\leq x$ and $y\leq 1/2$ as $\abs{C}\leq n/2$, an optimization of the above gives $e(X,Y)<dn^2$, which contradicts the definition of $d$.
\end{proof}

\begin{proof}[Proof of \Cref{lem:large_factor}.]
We construct a sequence of graphs $G=G_1\supseteq G_2\supseteq \dots\supseteq G_r$ in the following way. 
Let $d_i$ be the density of $G_i$ and $n_i$ be the order of $G_i$. 
We apply \Cref{lem:DensityIncrement} to $G_i$. 

If the first outcome occurs, then we let $G_{i+1}=G_i[A\cup B]$, for $A$ and $B$ as in \Cref{lem:DensityIncrement}.

If the second outcome occurs, then we stop the process here. 

Now note that we have constructed finitely many such graphs $G_i$, as $n_1,n_2,\dotsc,n_i,\dotsc$ is a decreasing sequence of positive integers.
Also note that for every $i$, we have immediately from \Cref{lem:DensityIncrement} that $d_in_i^{\delta} \geq d_{i-1}n_{i-1}^{\delta}$. 
Suppose that the process stopped after constructing $G_r$. Then
\begin{align*}
d_rn_r^2 &= d_r n_0^2 \prod_{i=0}^{r-1} \frac{n_{i+1}^2}{n_{i}^2}  \\
 &\geq d_r n_0^2 \Bigl(\prod_{i=0}^{r-1} \frac{d_{i}}{d_{i+1}}\Bigr)^{1/\delta}  \\
 &= d_r^{1-1/\delta} n_0^2 d_0^{1/\delta} \\
 &\geq n_0^2 d_0^{1/\delta}
\end{align*}

By definition of $r$ and \Cref{thm:criterion_k-factor}, the graph $G_r$ contains a $k$-factor, which has $kn_r/2=d_rn_r^2/200\geq d_0^{1/\delta}n_0^2/200$ edges, which finishes the proof.

\end{proof}

\section{On the maximal number of colours in an interval colouring of a planar graph}
\label{sec:main-proof}

First, we would like to emphasize that \Cref{thm:3/2result} is tight, as shown by Axenovich's construction \cite{Axenovich2002} of, for every natural integer $k$, a planar graph $G_k$ on $2k+2$ vertices satisfying $t(G_k) \geq 3k+1$.
At the end of this section, we extend this family to a larger collection of graphs, all of which demonstrate the tightness of this result.

We will now prove the following slightly more general result than \Cref{thm:3/2result}.

\begin{theorem}
\label{thm:generalisation-interval-colourable}
    Let $k > 1$ and $G$ be a graph on $n$ vertices such that every subgraph $H\subseteq G$ satisfies $|E(H)|\leq k(|V(H)|-2)$. 
    If $G$ admits an interval colouring, then $t(G)\leq (k/2)n + 1 - k$.
\end{theorem}

Note that a subgraph of a planar graph is also a planar graph, and that any planar graph $H$ satisfies $|E(H)|\leq 3(|V(H)|-2)$, so plugging $k=3$ into \Cref{thm:generalisation-interval-colourable} gives \Cref{thm:3/2result}. 
Moreover, note that \Cref{thm:generalisation-interval-colourable} is non-trivial only for $k \leq 4$ as otherwise it does not give any improvement over the result of Asratian and Kamalian \cite{asratian1994investigation} that $t(G)\leq 2n$ for any graph $G$ of order $n$.

We now prove \Cref{thm:generalisation-interval-colourable}.
The key result in our proof will be the following lemma.

\begin{lemma}
\label{lem:planar-key-lemma}
    Assume connected graph $G$ with interval colouring $c$ has an edge $e=vw$ such that $e$ is the unique edge to receive colour $c(e)=c_0$, and $c_0$ is neither the minimal nor maximal value taken by $c$.
    Then we may partition $V(G)=V_1\cup \{v,w\}\cup V_2$ such that every edge $f_1\in E(G)$ incident to some $x\in V_1$ satisfies $c(f_1)<c_0$, and every edge $f_2\in E(G)$ incident to some $x\in V_2$ satisfies $c(f_2)>c_0$.

    In particular, $e(V_1,V_2)=0$ and the induced colourings of $G[V_1\cup\{v,w\}]$ and $G[V_2\cup\{v,w\}]$ are interval.
\end{lemma}

\begin{proof}
    Firstly, as the colour $c_0$ of edge $e$ is unique, every vertex $x\in V(G)\setminus\{v,w\}$ has edges taking colours in an interval not including $c_0$, i.e. has colours either all less or all more than $c_0$.
    This naturally defines our parts $V_1,V_2$ to be those vertices whose edges have colours all less or all more than $c_0$ respectively.

    We also see that $e(V_1,V_2)=0$, as if we had some $f\in E(V_1,V_2)$, then we would have both that $c(f)<c_0$ and $c(f)>c_0$, a contradiction.

    To see that the induced colouring on $G[V_1\cup\{v,w\}]$ is interval, consider deleting every edge of $G$ with colour greater than $c_0$ to find a new graph $G'$.
    By the above results, $G[V_1\cup\{v,w\}]$ is a connected component of $G'$, and the induced colouring on $G'$ is interval.
    Thus the induced colouring on $G[V_1\cup\{v,w\}]$ is interval, and likewise for $G[V_2\cup\{v,w\}]$.
    This completes the proof.
\end{proof}

\begin{proof}[Proof of \Cref{thm:generalisation-interval-colourable}.]
Assume for contradiction that $G$ is a minimal counterexample, i.e. $n$ is minimal such that there is a graph $G$ on $n$ vertices with $t(G) > (k/2)n + 1 - k$.
We first claim that if such a graph $G$ is interval coloured with colours $1,\dotsc,t=t(G)$, then there must be some colour $c_0$ with $1 < c_0 < t$ such that there is a unique edge of $G$ receiving colour $c_0$.
Indeed, if this were not the case, then all colours except perhaps 1 and $t$ would occur at least twice, and thus $e(G)\geq 2(t-2)+2=2t-2$.
But by our assumption for contradiction, $e(G)\geq 2t-2 > k(n - 2) \geq e(G)$, which is impossible.

Therefore we may define $v,w,V_1$ and $V_2$ as in \Cref{lem:planar-key-lemma}, and let $G_1\defined G[V_1\cup\set{v,w}]$, and $G_2\defined G[V_2\cup\set{v,w}]$.
Then we know that $\abs{V_1}+\abs{V_2}=n+2$ and $t(G_1)+t(G_2) \geq t+1$.

By minimality of $G$, \Cref{thm:generalisation-interval-colourable} applies to $G_1$ and $G_2$. 
Then $t \leq t(G_1)+t(G_2)-1\leq (k/2)(\abs{V_1}+\abs{V_2})+1-2k=(k/2)n+1-k$, as required for a contradiction.
\end{proof}

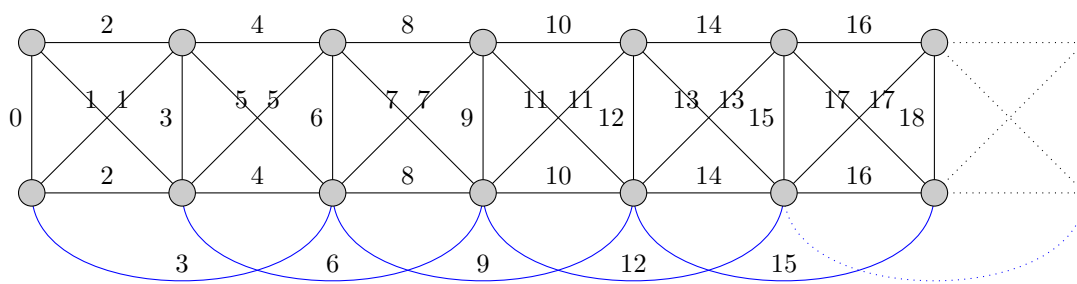
\begin{figure}[h]
\centering
\begin{tikzpicture}[scale=2,darkstyle/.style={circle,draw,fill=gray!40,minimum size=10}]
\path[-] 
    (1,1) edge ["0"] (1,2) 
        edge ["2"] (2,1) 
        edge ["1"] (2,2)
    (2,1) edge ["3"] (2,2) 
        edge ["4"] (3,1) 
        edge ["5"] (3,2)
    (3,1) edge ["6"] (3,2) 
        edge ["8"] (4,1) 
        edge ["7"] (4,2)
    (4,1) edge ["9"] (4,2) 
        edge ["10"] (5,1) 
        edge ["11"] (5,2)
    (5,1) edge ["12"] (5,2) 
        edge ["14"] (6,1) 
        edge ["13"] (6,2)
    (6,1) edge ["15"] (6,2) 
        edge ["16"] (7,1) 
        edge ["17"] (7,2)
    (7,1) edge ["18"] (7,2)
        edge [dotted] (8,1)
        edge [dotted] (8,2)
    (1,2) edge ["2"] (2,2) 
        edge ["1"] (2,1)
    (2,2) edge ["4"] (3,2) 
        edge ["5"] (3,1)
    (3,2) edge ["8"] (4,2) 
        edge ["7"] (4,1)
    (4,2) edge ["10"] (5,2) 
        edge ["11"] (5,1)
    (5,2) edge ["14"] (6,2) 
        edge ["13"] (6,1)
    (6,2) edge ["16"] (7,2) 
        edge ["17"] (7,1)
    (7,2) edge [dotted] (8,2)
        edge [dotted] (8,1)
    (1,1) edge [bend right=90, draw=blue, "3"] (3,1)
    (2,1) edge [bend right=90, draw=blue, "6"] (4,1)
    (3,1) edge [bend right=90, draw=blue, "9"] (5,1)
    (4,1) edge [bend right=90, draw=blue, "12"] (6,1)
    (5,1) edge [bend right=90, draw=blue, "15"] (7,1)
    (6,1) edge [bend right=90, draw=blue, dotted] (8,1);
\foreach \x in {1,...,7} 
    \foreach \y in {1,2} 
        {\node [darkstyle] (\x\y) at (\x,\y) {};}
\end{tikzpicture}
\caption{A planar graph on $2s$ vertices attaining the maximum value of $t(G)$, shown here for $s=7$. Any subset of the blue, curved edges can be removed to find another graph attaining this maximum.}
\label{fig:maximal-example-labelled}
\end{figure}

The bound in \Cref{thm:3/2result}, i.e. \Cref{thm:generalisation-interval-colourable} with $k=3$, is tight, and achieved by the family of graphs described as follows.
These graphs are found by starting with a graph on $2s$ vertices, like that shown in \Cref{fig:maximal-example-labelled}, and removing some subset of the edges which are drawn in blue and curved in that figure.
Note that after removing any such set of edges, the colouring shown is still an interval colouring, as the curved edges are either the minimal or maximal colour at every vertex they are incident to.

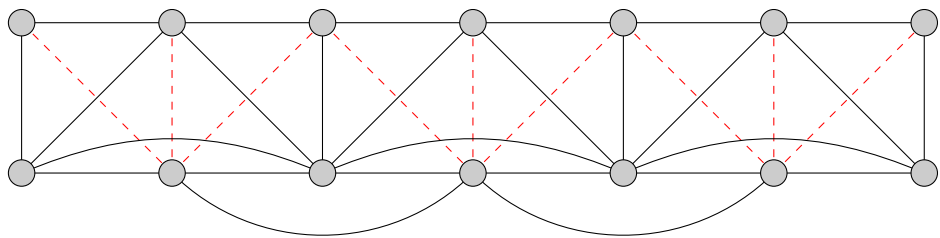
\begin{figure}[h]
\centering
\begin{tikzpicture}[scale=2,darkstyle/.style={circle,draw,fill=gray!40,minimum size=10}]
\path[-] 
    (1,1) edge (1,2) 
        edge (2,1) 
        edge (2,2)
    (2,1) edge [draw=red,dashed] (2,2) 
        edge (3,1) 
        edge [draw=red,dashed] (3,2)
    (3,1) edge  (3,2) 
        edge (4,1) 
        edge (4,2)
    (4,1) edge [draw=red,dashed] (4,2) 
        edge (5,1) 
        edge [draw=red,dashed] (5,2)
    (5,1) edge (5,2) 
        edge (6,1) 
        edge (6,2)
    (6,1) edge [draw=red,dashed] (6,2) 
        edge (7,1) 
        edge [draw=red,dashed] (7,2)
    (7,1) edge (7,2)
    (1,2) edge (2,2) 
        edge [draw=red,dashed] (2,1)
    (2,2) edge (3,2) 
        edge (3,1)
    (3,2) edge (4,2) 
        edge [draw=red,dashed] (4,1)
    (4,2) edge (5,2) 
        edge (5,1)
    (5,2) edge (6,2) 
        edge [draw=red,dashed] (6,1)
    (6,2) edge (7,2) 
        edge (7,1)
    (1,1) edge [bend left=23] (3,1)
    (2,1) edge [bend right=45] (4,1)
    (3,1) edge [bend left=23] (5,1)
    (4,1) edge [bend right=45] (6,1)
    (5,1) edge [bend left=23] (7,1);
\foreach \x in {1,...,7} 
    \foreach \y in {1,2} 
        {\node [darkstyle] (\x\y) at (\x,\y) {};}
\end{tikzpicture}
\caption{Demonstrating that our extremal examples are planar. The red dashed edges should loop around the left side of the graph.}
\label{fig:maximal-example-planar}
\end{figure}

To see that this graph is indeed planar, consider the drawing as shown in \Cref{fig:maximal-example-planar}. 
The dashed lines marked in red should be drawn to loop around the left-hand side of the graph.
As no red dashed line crosses any other such line, it can be seen that this drawing will be planar.

These examples all have $n$ even.
To see that the bound is still tight when $n$ is odd, take an extremal graph $G$ on $n-1$ vertices. 
Let $uv$ be an edge of $G$ given maximal colour $t$.
Add one new vertex, $w$, to $G$, and add a single edge of colour $t+1$ from $w$ to $v$.
This new graph then shows that the bound is tight (up to rounding) when $n$ is odd.

The extremal examples found by Axenovich \cite{Axenovich2002} correspond to the cases with all of the curved edges removed.

\section{Concluding remarks}

Borowiecka-Olszewska, Drgas-Burchardt, Javier-Nol and Zuazua \cite{borowiecka2021consecutive} defined an oriented graph to be \emph{consecutively colourable} if there exists a proper arc colouring such that for each vertex $v$, the colours of all out-arcs from $v$ and the colours of all in-arcs to $v$ form two intervals of integers. 
Moreover they conjectured that for every graph $G$, there exists an orientation of the edges that is consecutively colourable. It is easy to see that this conjecture implies that for every bipartite graph $G$, we have $\theta(G) \leq 2$. Therefore, \Cref{thm:LowerBoundTheta} also disproves Borowiecka-Olszewska et al's conjecture. \\

Moreover, we would like to point out that, given a constant $A \geq 1$, a similar lower bound than \Cref{thm:LowerBoundTheta} still holds (with a dependency of the constants in $A$) if instead of an interval-colouring, we ask for a proper colouring $c\from E\to \ZZ$ such that the set $C_x$ of colours of edges incident to vertex $x$ satisfies $\max(C_x)-\min(C_x) \leq Ad(x)$. \\

We would also like to point out that we have made no effort into optimizing the constants in this paper, therefore it is likely that one can improve the constant $8/9$ in \Cref{thm:UpperBoundTheta} by a more careful computation in the proof of \Cref{lem:DensityIncrement} for instance. \\

We proved \Cref{thm:UpperBoundTheta} by finding regular subgraphs with a large number of edges inside dense enough graphs in \Cref{lem:large_factor}. The problem of finding regular bipartite subgraphs with a large number of edges is interesting on its own, and it would be interesting to see if one could prove a stronger form of \Cref{lem:large_factor}. \\

In regards to the section on planar graphs, we suggest that these are in fact all examples which demonstrate that the bound in \Cref{thm:3/2result} is tight, and thus pose the following question.

\begin{Question}
\label{question:extremal-examples-3/2}
    Are all planar graphs $G$ on $n$ vertices for which $t(G)=3n/2-2$ of the form described above; that is, as shown in \Cref{fig:maximal-example-labelled} with some subset of the blue, curved edges removed?
\end{Question}

Another area of possible further work is to study whether \Cref{thm:generalisation-interval-colourable} is tight for values of $k$ other than $3/2$.

\begin{Question}
    For which values of $k$ is the bound in \Cref{thm:generalisation-interval-colourable} asymptotically tight? For such $k$, is there a description of all families of graphs for which the bound is (asymptotically) tight?
\end{Question}

\section*{Note added in proof} 

One day before the authors submitted this article to arXiv, Axenovich, Girão, Powierski, Savery and Tamitegama submitted a preprint \cite{https://doi.org/10.48550/arxiv.2303.04782} on arXiv where they independently proved that $\log(n)^{1/3-o(1)} \leq \theta(n) \leq n^{5/6+o(1)}$.

\bibliographystyle{abbrvnat}  
\renewcommand{\bibname}{Bibliography}
\bibliography{main}

\end{document}